\documentclass[12pt]{article}
\usepackage{mathrsfs}
\usepackage{amsfonts}

\begin{document}

\title{Global Existence of Solution for a Nonlinear Size-structured Population Model
with Distributed Delay in the Recruitment$^*$}
\author{Meng Bai$^{\dag}$\quad and \quad Shihe Xu \\
{\small School of Mathematics and Information Sciences, Zhaoqing
University,}\\
{\small Zhaoqing, Guangdong 526061 People's Republic of China. }}
\date{}
\maketitle

\begin{abstract}
In this paper we study a nonlinear size-structured population model
with distributed delay in the recruitment. The delayed problem is
reduced into an abstract initial value problem of an ordinary
differential equation in the Banach space by using the delay
semigroup techniques. The local existence and uniqueness of solution
as well as the continuous dependence on initial conditions are
obtained by using the general theory of quasi-linear evolution
equations in nonreflexive Banach spaces, while the global existence
of solution is obtained by the estimates of the solution and the
extension theorem.



\end{abstract}
\begin{figure}
\rule{50mm}{0.1mm}\\
$*$ \small This work is supported by NSF of China(11171295,11226182 and 11301474).\\
  $^{\dag}$ Corresponding author.
  \\Email addresses: baimeng.clare@aliyun.com(M.Bai), shihexu03@aliyun.com(S.Xu)
\end{figure}

\section{Introduction}\
In this paper, we study the following size-structured population
model with the distributed delay:
\begin{equation}
\left\{
\begin{array}{ll}
\begin{array}{rcl}
\displaystyle\frac{\partial n}{\partial t}+\frac{\partial
(\gamma(x,N[n(t)](x))n)}{\partial x}
&=&-\mu(x,N[n(t)](x))n+\displaystyle\int^{0}_{-\tau}
R[n(t+\sigma)](x)d\sigma\\[0.3cm]&&
\quad \mbox{for}\;\;x\in[0,\infty),t\geq0,\end{array}\\
[0.2cm] \displaystyle \gamma(0,N[n(t)](0))n(t,0)=0\quad
\mbox{for}\;\;t\geq0,\\
[0.2cm] \displaystyle n(\sigma,x)=\hat{n}(\sigma,x)\quad
\mbox{for}\;\;\sigma\in[-\tau,0],x\in[0,\infty).
\end{array}
\right.
\end{equation}
Here the unknown function $n(t,x)$ denotes the density of
individuals of size $x\in[0,\infty)$ at time $t\in[0,\infty)$,
$N[n(t)](x)$ is the {\em environment} or the {\em interaction
variable} (see\cite{AJ}) experienced by an individual of size $x$
when the population density is $n(t,y)$, $\gamma(x,N)$ and
$\mu(x,N)$ is the growth rate and the mortality rate of an
individual of size $x$ when the environment or the interaction
variable is $N$. We assume that there is a time lag $-\sigma$ in the
process of the recruitment, and $\sigma\in[-\tau,0]$, where $\tau>0$
is a constant denoting the maximal time lag. Mover precisely,
$R[n(t+\sigma)](x)$ is the portion of the recruitment of the new
individuals of size $x$ at time $t$ which come from the individuals
of size $y\in[0,\infty)$ at time $t+\sigma$,
$\int^{0}_{-\tau}R[n(t+\sigma)](x)$ is the entire recruitment of the
new individuals of size $x$ at time $t$. Besides,
$\hat{n}(\sigma,x)$ are a given function defined in $[-\tau,0]\times
[0,\infty)$. Later on we shall denote
\begin{equation}
\hat{n}_{0}(x)=\hat{n}(0,x)\;\;\; \hbox{for}\;\;x\in [0, \infty).
\end{equation}

The nonlinear size-structured population model without distributed
delay has been studied in \cite{AJ} and \cite{FGH}. Global existence
of solution has been obtained in \cite{AJ} and the asymptotic
behavior has been studied in \cite{FGH}. In the model (1.1), we
consider not only the environment but also the distributed delay in
the recruitment. More precisely, the distributed delay in the
recruitment is given by the time lag between conception and birth or
laying and hatching of the parasite eggs (see \cite{AOH}). Moreover,
unlike the non-distributed delay case, the time lag considered here
can change from $0$ to $\tau$, i.e., it is distributed in the
interval $[0,\tau]$. The linear age/size-structured population
models with delay in the birth process were studied in \cite{Pi},
\cite{PT} and \cite{MBSX}. Recently, some different nonlinear
age/size-structured population models with the distributed delay
were studied in \cite{XD} and \cite{ALAP}.

As an extension and development of the work in \cite{AJ}, we reduce
the distributed delayed problem (1.1) into an abstract initial value
problem of an ordinary differential equation in the Banach space and
obtain the global existence of the solution by using the delay
semigroup techniques, the general theory of quasi-linear evolution
equation in nonreflexive Banach spaces, the estimates of the
solution and the extension theorem.

For a given $K>0$, let $\Omega:=[0,\infty)\times[0,K]$. Then,
throughout this paper, $N$, $\mu(x,N)$, $\gamma(x,N)$, and $R$ are
supposed to satisfy the following conditions:

$(A.1)$\ \ Let $M_{n}:=\{v\in L^{\infty}(0,\infty),v'\in
W^{n,1}(0,\infty)\}$. Then, the operator $N:=L^{1}\rightarrow M_{0}$
and $N:=W^{1,1}\rightarrow M_{1}$, satisfies $N[0]=0$, and it is
Lipschitzian in the norms
$\|\cdot\|_{0}=\|\cdot\|_{\infty}+\|D(\cdot)\|_{L^{1}}$ and
$\|\cdot\|_{1}=\|\cdot\|_{\infty}+\|D(\cdot)\|_{W^{1,1}}$.

$(A.2)$\ \ $\mu(x,N)$ is a non-negative $\mathcal{C}^{1}-$ function,
and $\mu$, $\mu_{x}$ and $\mu_{N}$ are uniformly bounded by
$\mu^{0}$, $\mu^{0}_{x}$ and $\mu^{0}_{N}$, respectively, for all
$(x,N)\in\Omega$. Moreover, its partial derivatives $\mu_{x}$ and
$\mu_{N}$ are Lipschitzian functions with respect to $N$.

$(A.3)$\ \ $\gamma(x,N)$ is a strictly positive $\mathcal{C}^{2}-$
function for all $x,N\geq 0$, upper bounded by $\gamma^{0}>0$ and
lower bounded by $\gamma_{0}>0$ for all $(x,N)\in\Omega$. Moreover,
for all $(x,N)$, $|\gamma_{x}|$ and $|\gamma_{N}|$ are upper bounded
by $\gamma^{0}_{1}$, $|\gamma_{xx}|$, $|\gamma_{xN}|$ and
$|\gamma_{NN}|$ are upper bounded by $\gamma^{0}_{2}$. Finally,
$\gamma_{xx}(x,N)$, $\gamma_{xN}(x,N)$ and $\gamma_{NN}(x,N)$ are
Lipschitzian functions with respect to $N$.

$(A.4)$ The positive operator $R:L^{1}([-\tau,0],L^{1}(0,\infty))
\rightarrow L^{1}([-\tau,0], W^{1,1}(0,\infty))$ and
$R:L^{1}([-\tau,0], L^{1}(0,\infty))\rightarrow L^{1}([-\tau,0],
L^{\infty}(0,\infty))$ satisfies $R[0]=0$. For \\$\tilde{u}\in
L^{1}([-\tau,0],L^{1}(0,\infty))$,
$$\displaystyle\Big\|\int^{0}_{-\tau}R[\tilde{u}(\sigma)]d\sigma
\Big\|_{L^{1}(0,\infty)}\leq
R_{0}\int^{\infty}_{0}\int^{0}_{-\tau}\tilde{u}(\sigma,y)d\sigma
dy,$$
$$\displaystyle\Big\|\int^{0}_{-\tau}R[\tilde{u}(\sigma)]d\sigma
\Big\|_{W^{1,1}(0,\infty)}\leq
R_{1}\int^{\infty}_{0}\int^{0}_{-\tau}\tilde{u}(\sigma)d\sigma
dy,$$and
$$\Big\|\int^{0}_{-\tau}R[\tilde{u}(\sigma)]d\sigma
\Big\|_{\infty}\leq
R_{2}\int^{\infty}_{0}\int^{0}_{-\tau}\tilde{u}(\sigma,y)d\sigma
dy.$$ Moreover, for $\tilde{u}_{1},\tilde{u}_{2}\in E_{1},$ where
$E_{1}$ is an arbitrary bounded subset in
$L^{1}([-\tau,0],L^{1}(0,\infty))$,
$$\displaystyle\Big\|\int^{0}_{-\tau}R[\tilde{u}_{1}(\sigma)]d\sigma-
\int^{0}_{-\tau}R[\tilde{u}_{2}(\sigma)]d\sigma\Big\|_{L^{1}(0,\infty)}\leq
L_{R}\int^{\infty}_{0}\int^{0}_{-\tau}|\tilde{u}_{1}(\sigma,y)-\tilde{u}_{2}(\sigma,y)|d\sigma
dy,$$ and
$$\displaystyle\Big\|\int^{0}_{-\tau}R[\tilde{u}_{1}(\sigma)]d\sigma-
\int^{0}_{-\tau}R[\tilde{u}_{2}(\sigma)]d\sigma\Big\|_{W^{1,1}(0,\infty)}\leq
L_{R_{x}}\int^{\infty}_{0}\int^{0}_{-\tau}|\tilde{u}_{1}(\sigma,y)-\tilde{u}_{2}(\sigma,y)|d\sigma
dy.$$

The layout of the rest part is as follows. In Section 2 we reduce
the model (1.1) into an abstract initial value problem of an
ordinary differential equation in the Banach space. In Section 3 we
shall prove the local existence and uniqueness of solution of the
the model (1.1) by using the general theory of quasi-linear
evolution equations in nonreflexive Banach spaces. In Section 4 we
obtain the continuous dependence on initial conditions and the
positivity of solutions. In Section 5 we obtain the global existence
of solution by the estimates of the solution and the extension
theorem. In Section 6 we give the typical examples of the operators
$N$ and $R$.

\section{Reduction}
\setcounter{equation}{0}

In this section we reduce the problem $(1.1)$ into an abstract
Cauchy problem. We refer the reader to see $\cite{Pi}$, $\cite{PT}$
and $\cite{BP}$ for similar reductions.

First, we introduce the following Banach spaces:
\begin{eqnarray}
 &&X:=L^{1}(0,\infty),\quad \mbox{with norm}\;\;
  \|u\|_{X}=\int^{\infty}_{0}|u(x)|dx,\nonumber\\
&&Y:=\Big\{u\in W^{1,1}(0,\infty):
 u(0)=0\Big\},\quad \mbox{with norm}\;\;
  \|u\|_{Y}=\int^{\infty}_{0}|u(x)|dx+\int^{\infty}_{0}|u'(x)|dx,\nonumber\\
&&E:=L^{1}([-\tau,0],X),\quad \mbox{with norm}\;\;
  \|u\|_{E}=\displaystyle\int^{\infty}_{0}\int^{0}_{-\tau}|u(\sigma,x)|d\sigma dx.\nonumber\end{eqnarray}
For a given $v\in X$, let $A(v):Y\rightarrow X$ be the following
linear operator:
\begin{eqnarray}
&&A(v)u:=(\gamma(\cdot,N^{v})u)',\;\;\hbox{for}\;\;u\in Y,
\end{eqnarray}
where $N^{v}:=N[v]$. It is obvious that for a given $v\in X$,
$A(v)\in \mathcal{L}(Y,X)$. We denote by $f_{1}:X\rightarrow X$ and
$f_{2}:E\rightarrow X$, respectively, the following nonlinear
operators:
\begin{eqnarray}
 &&f_{1}(u):=-\mu(\cdot,N^{u})u,\;\;\hbox{for}\;\;u\in
 X,\nonumber\\
&&f_{2}(\tilde{u}):= \int^{0}_{-\tau}R[\tilde{u}(\sigma)]d\sigma,
\;\;\hbox{for}\;\;\tilde{u}\in E,\nonumber
\end{eqnarray}
where $N^{u}:=N[u]$.

Using these notations, we rewrite the model $(1.1)$ into the
following abstract initial value problem for a retarded differential
equation in the Banach space $X$:
\begin{equation}
\left\{
\begin{array}{ll}
\displaystyle\frac{d n(t)}{d
t}+A(n(t))n(t)=f_{1}(n(t))+f_{2}(n_{t}), &\;\;t\geq0\\
[0.2cm] \displaystyle n(0)=\hat{n}_{0},
   \\[0.2cm]
\displaystyle n_{0}=\hat{n},
\end{array}
\right.
\end{equation}
where $n:[0,+\infty)\rightarrow X $ is defined as $n(t):
=n(t,\cdot)$ and $n_{t}:[-\tau,0]\rightarrow X $ is defined as
$n_t(\sigma):=n(t+\sigma),\sigma\in[-\tau,0].$

Next, we introduce the following operators in the Banach space $E$:
\begin{eqnarray}
&&(G\tilde{u})(\sigma):=-\frac{d}{d\sigma}\tilde{u},\quad \hbox{with
domain}\;\;
D(G)=W^{1,1}([-\tau,0],X),\nonumber\\
&&Q \tilde{u}:=\tilde{u}(0),\quad\hbox{for}\;\;\tilde{u}\in
D(G).\nonumber
\end{eqnarray}
We note that $G\in\mathcal{L}(D(G),E)$ and
$Q\in\mathcal{L}(D(G),X)$. We now let
$$\mathbb{X}:=E\times X,\quad \mbox{with norm}\;\;
  \|(\tilde{u},u)\|_{\mathbb{X}}=\|\tilde{u}\|_{E}+\|u\|_{X}$$
and $$\mathbb{Y}:=\Big\{(\tilde{u},u)\in W^{1,1}([-\tau,0],X)\times
W^{1,1}(0,\infty), Q \tilde{u}=u ,u(0)=0\Big\},$$ $$\quad \mbox{with
norm}\;\;\|(\tilde{u},u)\|_{\mathbb{Y}}=\|\tilde{u}\|_{E}+\|\tilde{u}'_{1}\|_{E}+\|u\|_{Y},
$$
where $\tilde{u}'_{1}(\sigma,x)=\displaystyle\frac{\partial
\tilde{u}(\sigma,x)}{\partial\sigma} $. For a given
$\mathbf{V}=(\tilde{v},v)\in \mathbb{X}$, let
$\mathbf{A}(\mathbf{V}):\mathbb{Y}\rightarrow \mathbb{X}$ be the
following operator:
\begin{eqnarray}
&&\mathbf{A}(\mathbf{V})\mathbf{U}:=\left(
\begin{array}{cc}
G\;\;\;&0\\
0\;\;\;&A(v)
\end{array}
\right)\left(
\begin{array}{c}
\tilde{u}\\
u
\end{array}
\right),\;\;\mbox{for}\;\; \mathbf{U}=\left(
\begin{array}{c}
\tilde{u}\\
u
\end{array}
\right)\in\mathbb{Y}.
\end{eqnarray}
It is obvious that for given $\mathbf{V}=(\tilde{v},v)\in
\mathbb{X}$, $\mathbf{A}(\mathbf{V})\in \mathcal{L}( \mathbb{Y},
\mathbb{X})$. We also denote by $\mathbf{F}:\mathbb{X}\rightarrow
\mathbb{X}$, the following nonlinear operator:
\begin{eqnarray}
\mathbf{F}(\mathbf{U})=\left(
\begin{array}{c}
0\\
f_{1}(u)+f_{2}(\tilde{u})
\end{array}
\right),\quad \mbox{for}\;\; \mathbf{U}=\left(
\begin{array}{c}
\tilde{u}\\
u
\end{array}
\right)\in\mathbb{X}.
\end{eqnarray}

Using these notations, we see that the problem (2.1) can be
equivalently rewrite into the following abstract initial value
problem of an ordinary differential equation in the Banach space
$\mathbb{X}$:
\begin{equation}
\left\{
\begin{array}{ll}
\displaystyle \mathbf{U}'(t)+\mathbf{A}(\mathbf{U}(t))\mathbf{U}(t)
=\mathbf{F}(\mathbf{U}(t)), &\;\;t>0,\\
[0.2cm] \displaystyle  \mathbf{U}(0)=\mathbf{U}_{0},\\ [0.2cm]
\end{array}
\right.
\end{equation}
where $\mathbf{U}(t)=\left(
\begin{array}{c}
n_{t}\\
n(t)
\end{array}
\right)$ and $\mathbf{U}_{0}=\left(
\begin{array}{c}
\hat{n}\\
\hat{n}_{0}
\end{array}
\right)$, where $\hat{n}_{0}$ is the function defined in $(1.2)$.

To describe the relationship between the problem (2.2) and the
problem (2.5), we write down the following preliminary result:

\medskip
{\bf Lemma 2.1}\ \ {\em If for any initial condition
$(\hat{n},\hat{n}_{0})\in W^{1,1}([-\tau,0],X)\times Y$, there
exists a time $T>0$ such that the problem $(2.2)$ has a unique
solution $n\in C([-\tau,T], W^{1,1}(0,\infty))\cap C([0,T],Y)\cap
C^{1}([0,T],X)$, then the problem $(2.5)$ has a unique solution
$\mathbf{U}(t)=\left(
\begin{array}{c}
n_{t}\\
n(t)
\end{array}
\right)$ and $\mathbf{U}\in C([0,T],\mathbb{Y})\cap
C^{1}([0,T],\mathbb{X})$. Conversely, if for any initial condition
$\mathbf{U}_{0}\in \mathbb{Y}$, there exists a time $T>0$ such that
$(2.5)$ has a unique solution $\mathbf{U}\in C([0,T],\mathbb{Y})\cap
C^{1}([0,T],\mathbb{X})$, then $\mathbf{U}$ has the form
$\mathbf{U}(t)=\left(
\begin{array}{c}
n_{t}\\
n(t)
\end{array}
\right)$ for all $t\in[0,T]$, and by extending its second component
$n=n(t)$ to $[-\tau,T]$ such that $n(t)=\hat{n}$ for $t\in
[-\tau,0)$, we have that a unique solution of the problem $(2.2)$.}
\medskip

{\em Proof}:\ \  We only need to prove that for any initial
condition $\mathbf{U}_{0}\in \mathbb{Y}$, there exists a time $T>0$
such that the problem (2.5) has a unique solution $\mathbf{U}\in
C([0,T],\mathbb{Y})\cap C^{1}([0,T],\mathbb{X})$, then $\mathbf{U}$
has the form $\mathbf{U}(t)=\left(
\begin{array}{c}
n_{t}\\
n(t)
\end{array}
\right)$ for all $t\in[0,T]$. To see this we assume that
$\mathbf{U}(t)=\left(
\begin{array}{c}
U(t)\\
n(t)
\end{array}
\right)$ (for $t\in[0,T]$). Since $\mathbf{U}\in
C([0,T],\mathbb{Y})$, we have $QU(t)=n(t)$ for all $t\in[0,T]$,
i.e., $U(t)(0)=n(t)$ for all $t\in[0,T]$ (recall that $U(t)\in
W^{1,1}([-\tau,0],X)$ for every $t\in[0,T]$). Let
$v(t,\sigma)=U(t)(\sigma)$ (for $t\in[0,T]$ and $\sigma\in
[-\tau,0]$). Then from the equation satisfied by $U$ we see that $v$
satisfies the equation $\partial v/\partial t-\partial
v/\partial\sigma=0$, so that it is a leftward traveling wave, i.e.,
it has the form $v(t,\sigma)=g(t+\sigma)$ for some function $g=g(s)$
defined for $s\in[-\tau,T]$. For $s\in[0,T]$ we have $g(s)=v(s,0)
=U(s)(0)=n(s)$. Hence, for $t+\sigma\in[0,T]$ we have
$$
  U(t)(\sigma)=v(t,\sigma)=g(t+\sigma)=n(t+\sigma)=n_t(\sigma).
$$
Moreover, since $U(0)(\sigma)=\hat{n}(\sigma,\cdot)$ and
$U(t)(\sigma)=v(t,\sigma)= g(t+\sigma)$ for $-\tau\leq t+\sigma<0$,
we see that $g(s)=\hat{n}(s,\cdot)$ for $s\in [-\tau,0)$ and
$U(t)(\sigma)=\hat{n}(t+\sigma,\cdot)$ for $t+\sigma\in [-\tau,0)$.
Hence, by defining $n(t)=\hat{n}(t,\cdot)$ for $t\in [-\tau,0)$, we
see that $U(t)(\sigma)=n_t(\sigma)$ also holds for all $t\in[0,T]$
and $\sigma\in [-\tau,0]$. This proves the desired assertion.
$\quad\Box$
\medskip

\section{Local existence and uniqueness of
solution} \setcounter{equation}{0}

In this section we shall prove the local existence and uniqueness of
solution for the problem (2.5) by using the general theory of
quasi-linear evolution equations in nonreflexive Banach
spaces(see\cite{AJ} and \cite{KN}). For this purpose, we shall
verify the hypotheses in the Theorems which are proposed and proven
by Kobayasi and Sanekata(see \cite{KN} and Theorem I and Theorem II
of \cite{AJ}), and apply these theorems to the problem (2.5). In the
sequel we denote by $\mathbb{W}$ the open subset of $\mathbb{Y}$
which contained in the closed ball in $\mathbb{Y}$ with center $0$
and radius $r>\|\mathbf{U}_{0}\|_{\mathbb{Y}}$,
$S(\mathbb{X},M,\alpha)$ the set of all stable families of negative
generators of $C^{0}-$ semigroups in $\mathbb{X}$ with stability
index $(M,\alpha)$, and $\mathcal{B}(X,Y)$ the set of all bounded
linear operators for the real Banach space $X$ to the real Banach
space $Y$($\mathcal{B}(X):=\mathcal{B}(X,X)$).

\medskip
{\bf Lemma 3.1}\ \ {\em The space $\mathbb{Y}$ is densely and
continuously embedded in $\mathbb{X}$. There is an isomorphism
$\mathbf{S}$ of $\mathbb{Y}$ onto $\mathbb{X}$. Moreover, there
exists two positive constants $c_{1}$ and $c_{2}$ such that
$$
c_{1}\|\mathbf{U}\|_{\mathbb{Y}}\leq
\|\mathbf{S}\mathbf{U}\|_{\mathbb{X}}\leq
c_{2}\|\mathbf{U}\|_{\mathbb{Y}}.
$$}
\medskip
{\em Proof}:\ \ Since there exist a generator with domain
$\mathbb{Y}$ which generates a strongly continuous semigroup on
$\mathbb{X}$ (see Proposition 4.4 of \cite{PT}), by Corollary 2.5 in
Chapter 1 of \cite{AP}, the first assertion follows.

The isomorphism $\mathbf{S}: \mathbb{Y}\rightarrow \mathbb{X}$ is
$$ \mathbf{S}\left(
\begin{array}{c}
\tilde{u}\\
u
\end{array}
\right)=\left(
\begin{array}{c}
-\tilde{u}'_{1}+\tilde{u}\\
u'+u
\end{array}
\right),\quad \mbox{for}\;\; \left(
\begin{array}{c}
\tilde{u}\\
u
\end{array}
\right)\in\mathbb{Y}.
$$
$\mathbf{S}$ is bijective because for every
$\mathbf{F}=(\tilde{f},f)\in\mathbb{X}$, the equations
\begin{equation}
\left\{
\begin{array}{ll}
\displaystyle
-\frac{\partial}{\partial\sigma}\tilde{u}(\sigma,x)+\tilde{u}(\sigma,x)
=\tilde{f}(\sigma,x),\;\;-\tau<\sigma< 0,0<x<\infty,\\
[0.2cm] \displaystyle
\frac{du(x)}{dx}+u(x)=f(x),\;\;0<x< \infty,\\[0.2cm]
\displaystyle \tilde{u}(0,x)=u(x),\;\;0<x< \infty,\\[0.2cm]
u(0)=0,
\end{array}
\right.
\end{equation}
have a unique solution, given by
\begin{equation}
\left\{
\begin{array}{ll}
\displaystyle
\tilde{u}(\sigma,x)=e^{\sigma}u(x)+e^{\sigma}\int^{0}_{\sigma}e^{-\xi}\tilde{f}(\xi,x)d\xi,\\
[0.2cm] \displaystyle u(x)=e^{-x}\int^{x}_{0}e^{s}f(s)ds,
\end{array}
\right.
\end{equation}
and the unique solution $(\tilde{u},u)\in \mathbb{Y}$.

Let $\mathbf{U}=(\tilde{u},u)\in \mathbb{Y}$. It is easy to see that
$\|\mathbf{U}\|_{\mathbb{Y}}\geq\|\mathbf{S}\mathbf{U}\|_{\mathbb{X}}$.
This implies that $c_{2}=1$. Since
$$
\int^{\infty}_{0}|u(x)|dx\leq\int^{\infty}_{0}e^{-x}\int^{x}_{0}e^{s}|f(s)|dsdx
\leq
\int^{\infty}_{0}\Big(\int^{\infty}_{s}e^{-x}dx\Big)e^{s}|f(s)|ds
=\int^{\infty}_{0}|f(x)|dx
$$
and
\begin{eqnarray}
&&
\displaystyle\int^{\infty}_{0}\int^{0}_{-\tau}|\tilde{u}(\sigma,x)|d\sigma
dx \nonumber
\\&\leq& \int^{\infty}_{0}\int^{0}_{-\tau}|e^{\sigma}u(x)|d\sigma dx+
\int^{\infty}_{0}\int^{0}_{-\tau}\Big|e^{\sigma}\int^{0}_{\sigma}e^{-\xi}\tilde{f}(\xi,x)d\xi\Big|d\sigma
dx \nonumber
\\&\leq& \tau\int^{\infty}_{0}|u(x)|dx
+\int^{\infty}_{0}\int^{0}_{-\tau}\int^{0}_{\sigma}e^{\sigma-\xi}|\tilde{f}(\xi,x)|d\xi
d\sigma dx \nonumber
\\&\leq&\displaystyle \tau\int^{\infty}_{0}|f(x)|dx
+\tau\int^{\infty}_{0}\int^{0}_{-\tau}|\tilde{f}(\sigma,x)|d\sigma
dx,\nonumber
\end{eqnarray}
we have
$$
\begin{array}{rcl}
\displaystyle\int^{\infty}_{0}|u'(x)|dx
&=&\displaystyle\int^{\infty}_{0}|f(x)-u(x)|dx\\[0.1cm]
&\leq& \displaystyle
\int^{\infty}_{0}|u(x)|dx+\int^{\infty}_{0}|f(x)|dx \\[0.1cm]
&\leq& 2\int^{\infty}_{0}|f(x)|dx
\end{array}
$$
and
$$
\begin{array}{rcl}
\displaystyle\int^{\infty}_{0}\int^{0}_{-\tau}\Big|\frac{\partial
\tilde{u}(\sigma,x)}{\partial \sigma}\Big|d\sigma
dx&=&\displaystyle\int^{\infty}_{0}\int^{0}_{-\tau}|\tilde{u}(\sigma,x)-
\tilde{f}(\sigma,x)|d\sigma
dx\\[0.1cm]
&\leq& \displaystyle
\int^{\infty}_{0}\int^{0}_{-\tau}|\tilde{u}(\sigma,x)|d\sigma
dx+\int^{\infty}_{0}\int^{0}_{-\tau}|\tilde{f}(\sigma,x)|d\sigma
dx \\[0.1cm]
&\leq&\displaystyle\tau\int^{\infty}_{0}|f(x)|dx+(\tau+1)\int^{\infty}_{0}\int^{0}_{-\tau}
|\tilde{f}(\sigma,x)|d\sigma
dx.
\end{array}
$$
Then
$$
\begin{array}{rcl}
\displaystyle\|\mathbf{U}\|_{\mathbb{Y}}&=&\displaystyle
\int^{\infty}_{0}\int^{0}_{-\tau}|\tilde{u}(\sigma,x)|d\sigma
dx+\int^{\infty}_{0}\int^{0}_{-\tau}\Big|\frac{\partial
\tilde{u}(\sigma,x)}{\partial \sigma}\Big|d\sigma dx+
\int^{\infty}_{0}|u(x)|dx+\int^{\infty}_{0}|u'(x)|dx\\[0.1cm]
&\leq& \displaystyle
(2\tau+1)\int^{\infty}_{0}\int^{0}_{-\tau}|\tilde{f}(\sigma,x)|d\sigma
dx+(2\tau+3)\int^{\infty}_{0}|f(x)|dx\\[0.1cm]
&\leq& \displaystyle(2\tau+3)\|\mathbf{S}\mathbf{U}\|_{\mathbb{X}}.
\end{array}
$$
This implies that $\displaystyle c_{2}=\frac{1}{2\tau+3}$. This
completes the proof.$\quad\Box$

\medskip
{\bf Lemma 3.2}\ \ {\em For each $\mathbf{W}\in \mathbb{W}$,
$\mathbf{A}(\mathbf{W})$ is a linear operator in $\mathbb{W}$. Let
$T$ be a positive constant. For each $\rho\geq0$, there exist two
constants $M\geq1$ and $\alpha\geq0$ such that
$$(\mathbf{A}(\mathbf{V}(t))_{t\in[0,T]}\in
S(\mathbb{X},M,\alpha),\;\;\forall
\mathbf{V}(t)\in\mathcal{D}_{\rho},
$$
where $\mathcal{D}_{\rho}:=\{\mathbf{V}(t)\in
C([0,T],\mathbb{W}):\|\mathbf{V}(t)-\mathbf{V}(s)\|_
{\mathbb{X}}\leq\rho|t-s|,0\leq
s<t\leq T\}$.}
\medskip

{\em Proof}:\ \ From (2.1) and (2.3), the first assertion follows.
We denote by $\mathbf{I}=\left(
\begin{array}{cc}
\tilde{I}\;\;\;&0\\
0\;\;\;&I
\end{array}
\right)$ the identity operator in $\mathbb{X}$, where $\tilde{I}$
and $I$ represent the identity operators in $E$ and $X$,
respectively. For each $\mathbf{W}\in \mathbb{W}$, we introduce the
operator
$\mathbf{A}_{1}(\mathbf{W}):=\mathbf{A}(\mathbf{W})+\mathbf{I}$.
Then for each $\mathbf{W} \in \mathbb{W}$,
$\mathbf{A}(\mathbf{W})=\mathbf{A}_{1}(\mathbf{W})-\mathbf{I}$.
Since $\|\mathbf{I}\|_{\mathcal{L}(\mathbb{X})}\leq 1$, by Theorem
2.3 in Chapter 5 of \cite{AP}, the second assertion follows if we
prove the stability of the family
$(\mathbf{A}_{1}(\mathbf{V}(t))_{t\in[0,T]}$. Since for a given
$\mathbf{W}\in \mathbb{W}$, the operator
$\mathbf{A}_{1}(\mathbf{W})$ does not depend on $t$ we can take
$T=\infty$ and, in order to prove the stability of the family
$(\mathbf{A}_{1}(\mathbf{V}(t))_{t\in[0,T]}$ we only need to show
that for each $\mathbf{W}=(\tilde{w},w)\in \mathbb{W}$,
$-\mathbf{A}_{1}(\mathbf{W})$ generates a contraction semigroup. To
this end we prove that for each $\mathbf{W}=(\tilde{w},w)\in
\mathbb{W}$, $-\mathbf{A}_{1}(\mathbf{W})$ satisfy the conditions of
Hille-Yosida Theorem (see Theorem 3.1 in Chapter 1 of \cite{AP}).
Since for each $\mathbf{W}=(\tilde{w},w)\in \mathbb{W}$,
$-\mathbf{A}_{1}(\mathbf{W})$ is closed,
$\overline{D(-\mathbf{A}_{1}(\mathbf{W}))}=\overline{\mathbb{Y}}=\mathbb{X}$
and the resolvent set $\rho(-\mathbf{A}_{1}(\mathbf{W}))$ of
$-\mathbf{A}_{1}(\mathbf{W})$ contains $\mathbb{R}^{+}$ (see
Proposition 2.1 of \cite{G} and Proposition 3.2 of \cite{PT}), we
only need to show that for every $\lambda>0$,
$\displaystyle\|R(\lambda,-\mathbf{A}_{1}(\mathbf{W}))
\|\leq\frac{1}{\lambda}$.
For $\mathbf{F}\in \mathbb{X}$, let
$U=R(\lambda,-\mathbf{A}_{1}(\mathbf{W}))\mathbf{F}$. Then
$\mathbf{U}$ satisfies the equation
\begin{eqnarray}
(\lambda
\mathbf{I}+\mathbf{A}_{1}(\mathbf{W}))\mathbf{U}=\mathbf{F}.
\end{eqnarray}
By writing $\mathbf{U}=(\tilde{u}(\sigma,x),u(x))$ and
$\mathbf{F}=(\tilde{f}(\sigma,x),f(x))$, we see that the above
equation can be rewritten as follows:
\begin{equation}
\left\{
\begin{array}{ll}
\displaystyle\lambda
\tilde{u}(\sigma,x)+\tilde{u}(\sigma,x)-\frac{\partial}{\partial\sigma}\tilde{u}(\sigma,x)
=\tilde{f}(\sigma,x),\;\;-\tau<\sigma< 0,0<x<m,\\
[0.2cm] \displaystyle
\lambda u(x)+u(x)+\frac{d}{dx}(\gamma(x,N^{w})u(x))=f(x),\;\;0<x< m,\\[0.2cm]
\displaystyle U(0,x)=u(x),\;\;0<x< m,\\[0.2cm]
u(0)=0.\\[0.2cm]
\end{array}
\right.
\end{equation}Then we
have that
\begin{equation}
\tilde{u}(\sigma,x)=e^{(\lambda+1)\sigma}u(x)+e^{(\lambda+1)\sigma}
\int^{0}_{\sigma}e^{-(\lambda+1)\xi}\tilde{f}(\xi,x)d\xi,
\end{equation}
and
\begin{equation} \displaystyle u(x)= \displaystyle
E_{\lambda}(x)\int^{x}_{0}(E_{\lambda}(s)\gamma(s,N^{w}))^{-1}f(s)ds,\nonumber
\end{equation}
where $
E_{\lambda}(x)=\displaystyle\exp\left\{-\int^{x}_{0}\frac{\lambda+
1+\gamma'(s,N^{w})}{\gamma(s,N^{w})}ds\right\}.
$ We deduce an useful expression of
$R(\lambda,-\mathbf{A}_{1}(\mathbf{W}))$. From (3.4), (3.5) and
(3.6), we have that
\begin{equation}
R(\lambda,-\mathbf{A}_{1}(\mathbf{W}))=\left(
\begin{array}{cc}
R(\lambda,-G_{0})\;\;\;&\varepsilon_{\lambda}R(\lambda,-A_{1}(w))\\
0\;\;\;&R(\lambda,-A_{1}(w))
\end{array}
\right),
\end{equation}
where $\varepsilon_{\lambda}:=e^{(\lambda+1)\sigma}$ for
$\sigma\in[-\tau,0)$, $A_{1}(w):=A(w)+I$ and $G_{0}$ is the
following operator in the Banach space $E=L^{1}([-\tau,0],X)$:
\begin{eqnarray}
&&(G_{0}\tilde{u})(\sigma):=-\frac{d}{d\sigma}\tilde{u}+\tilde{u},\quad
\hbox{with domain}\;\; D(G_{0})=\Big\{\tilde{u}\in
W^{1,1}([-\tau,0],X):
 \tilde{u}(0,x)=0\Big\}.\nonumber
\end{eqnarray}
Since $\mathbf{W}=(\tilde{w},w)\in \mathbb{W}$, there exists a
bounded set $W$ contained in $Y$ such that $w\in W$. By the proof of
H2 in \cite{AJ}, we have that for each $w\in W$, $-A(w)$ generates a
contraction semigroup $\{T(t)\}_{t\geq 0}$. Since
$-A_{1}(w)=-A(w)-I$, by the theory of the rescaled semigroups(see
Example II 2.2 of \cite{EN}), we have that for each $w\in W$,
$-A_{1}(w)$ generates a semigroup $\{S(t)\}_{t\geq 0}$ such that
$S(t)=e^{-t}T(t)$ for $t\geq 0$. Then $\|S(t)\|\leq e^{-t}$ for
$t\geq 0$. By Corollary 3.8 in Chapter 1 of \cite{AP}, we have that
for every $\lambda>0$,
$\displaystyle\|R(\lambda,-A_{1}(w))\|_{\mathcal{L}(X)}\leq
\frac{1}{\lambda+1}$. Hence, for
$\mathbf{F}=(\tilde{f}(\sigma,x),f(x))\in \mathbb{X}$ and
$\lambda>0$,
\begin{eqnarray}
&& \displaystyle
\|R(\lambda,-\mathbf{A}_{1}(\mathbf{W}))\mathbf{F}\|_{\mathbb{X}}\nonumber
\\&=&
\|R(\lambda,G_{0})F+\varepsilon_{\lambda}R(\lambda,-A_{1}(w))f\|_{E}+
\|R(\lambda,-A_{1}(w))f\|_{X}\nonumber\\
&\leq&
\displaystyle\int^{\infty}_{0}\int^{0}_{-\tau}e^{(\lambda+1)\sigma}
\int^{0}_{\sigma}e^{-(\lambda+1)\xi}|\tilde{f}(\xi,x)|d\xi
d\sigma
dx\nonumber\\& &+\int^{\infty}_{0}\int^{0}_{-\tau}e^{(\lambda+1)\sigma}
|R(\lambda,-A_{1}(w))
f(x)|d\sigma dx+\frac{1}{\lambda+1}\|f\|_{X}\nonumber\\
&\leq&
\displaystyle\int^{\infty}_{0}\int^{0}_{-\tau}e^{-(\lambda+1)\xi}|\tilde{f}(\xi,x)|
\int^{\xi}_{-\tau}e^{(\lambda+1)\sigma}d\sigma d\xi dx\nonumber\\&
&+
\int^{0}_{-\tau}e^{(\lambda+1)\sigma}d\sigma\int^{\infty}_{0}|R(\lambda,
-A_{1}(w))f(x)|dx+
\frac{1}{\lambda+1}\|f\|_{X}\nonumber\\
&\leq&\displaystyle\frac{1}{\lambda+1}\|\tilde{f}\|_{E}
+\frac{1}{(\lambda+1)^{2}}\|f\|_{X}+\frac{1}{\lambda+1}\|f\|_{X}\nonumber\\
&\leq&\displaystyle\frac{1}{\lambda+1}\|\tilde{f}\|_{E}
+(\frac{1}{(\lambda+1)^{2}}+\frac{1}{\lambda+1})\|f\|_{X}\nonumber\\
&\leq&\displaystyle\frac{1}{\lambda}\|\mathbf{F}\|_{\mathbb{X}}.\nonumber\
\end{eqnarray}
This completes the proof. $\quad\Box$

\medskip
 {\bf Lemma 3.3}\ \ {\em
For any $\mathbf{W}\in \mathbb{W}$, there exists an operator
$\mathbf{B}(\mathbf{W})$ such that
$$
\mathbf{S}\mathbf{A}(\mathbf{W})\mathbf{S}^{-1}=\mathbf{A}
(\mathbf{W})+\mathbf{B}(\mathbf{W}),
\mathbf{W}\in \mathbb{W},
$$
where $\mathbf{S}$ is the isomorphism defined in Lemma 3.1.
Moreover, there exist two positive numbers $\lambda_{\mathbf{B}}$
and $\mu_{\mathbf{B}}$ such that
\begin{eqnarray}
\|\mathbf{B}(\mathbf{W})\|_{\mathbb{X}}\leq
\lambda_{\mathbf{B}},\;\;\hbox{for}\;\;\mathbf{W}\in \mathbb{W},
\end{eqnarray}
and
\begin{eqnarray}
\|\mathbf{B}(\mathbf{W}_{1})-\mathbf{B}(\mathbf{W}_{2})\|_{\mathbb{X}}
\leq
\mu_{\mathbf{B}}\|\mathbf{W}_{1}-\mathbf{W}_{2}\|_{\mathbb{Y}}
,\;\;\hbox{for}\;\;\mathbf{W}_{1},\mathbf{W}_{2}\in
\mathbb{W}.
\end{eqnarray}}
\medskip
{\em Proof}:\ \ For a given $\mathbf{W}=(\tilde{w},w)\in
\mathbb{W}$, we have that
$$
\mathbf{B}(\mathbf{W})\mathbf{U}=(\mathbf{S}
\mathbf{A}(\mathbf{W})-\mathbf{A}(\mathbf{W})\mathbf{S})
\mathbf{S}^{-1}\mathbf{U},\;\;\hbox{for}\;\;\mathbf{U}\in
\mathbb{X}.
$$
To find the concrete expression of $\mathbf{B}(\mathbf{W})$, we
compute the "commutator"
$\mathbf{S}\mathbf{A}(\mathbf{W})-\mathbf{A}(\mathbf{W})\mathbf{S}$:
$$
(\mathbf{S}\mathbf{A}(\mathbf{W})-\mathbf{A}(\mathbf{W})\mathbf{S})
\mathbf{U}=\left(
\begin{array}{c}
0\\
D^{2}(\gamma(\cdot,N^{w}))u+D(\gamma(\cdot,N^{w}))u'
\end{array}
\right),\;\;\hbox{for}\;\;\mathbf{U}=\left(
\begin{array}{c}
\tilde{u}\\
u
\end{array}
\right)\in \mathbb{Y}
$$
where $\displaystyle
D(\gamma(x,N^{w}))=\gamma'_{1}(x,N^{w})+\gamma'_{2}(x,N^{w})\frac{d
N^{w}}{dx}$.

Let $\mathbf{U}\in \mathbb{X}$ and
$\mathbf{V}=\mathbf{S}^{-1}\mathbf{U}\in \mathbb{Y}$. Then,
$$
\mathbf{B}(\mathbf{W})\mathbf{U}=
(\mathbf{S}\mathbf{A}(\mathbf{W})-\mathbf{A}(\mathbf{W})\mathbf{S})
\mathbf{S}^{-1}\mathbf{U}=
(\mathbf{S}\mathbf{A}(\mathbf{W})-\mathbf{A}(\mathbf{W})\mathbf{S})
\mathbf{V}
$$
The first component of $\mathbf{B}(\mathbf{W})$ is zero. The second
component of $\mathbf{B}(\mathbf{W})$ is similar as the operator
$B(w)$ in the proof of H3 in \cite{AJ}. Under (A.1) and (A.3), we
can obtain (3.8) and (3.9) by using the same method. This completes
the proof.$\quad\Box$

{\em Remark 3.1:}\ \  Since for a given $\mathbf{W}\in \mathbb{W}$,
the operator $\mathbf{B}(\mathbf{W})$ does not depend on $t$,
$\mathbf{B}(\mathbf{W}):[0,T_{0}]\rightarrow\mathcal{B}(\mathbb{X})$
is a strongly measurable operator valued function with
$T_{0}=\infty$.

\medskip
{\bf Lemma 3.4}\ \ {\em For each $\mathbf{W}\in \mathbb{W}$,
$D(\mathbf{A}(\mathbf{W}))\supset \mathbb{Y}$ and
$\mathbf{A}(\mathbf{W})\in \mathcal{B}(\mathbb{Y},\mathbb{X})$.
Moreover, there exists a constant $\mu_{\mathbf{A}}$ such that
\begin{eqnarray}
\|\mathbf{A}(\mathbf{W}_{1})-\mathbf{A}(\mathbf{W}_{2})\|_
{\mathbb{Y},\mathbb{X}}
\leq
\mu_{\mathbf{A}}\|\mathbf{W}_{1}-\mathbf{W}_{2}\|_
{\mathbb{X}},\;\;\hbox{for}
\;\;\mathbf{W}_{1},\mathbf{W}_{2}\in
\mathbb{W}.
\end{eqnarray}}
\medskip
{\em Proof}:\ \ For any $\mathbf{W}=(\tilde{w},w)\in\mathbb{W}$,
$D(\mathbf{A}(\mathbf{W}))= \mathbb{Y}$. Let $\mathbf{U}\in
\mathbb{Y}$, we have that
$$
\begin{array}{rcl}
\displaystyle\|\mathbf{A}(\mathbf{W})\mathbf{U}\|_{\mathbb{X}}
&=&\displaystyle
\Big\|G\tilde{u}\Big\|_{E}+\|(A(w))u\|_{X}\\[0.1cm]&\leq&
\displaystyle\Big\|\frac{\partial\tilde{u}(\sigma,x)}{\partial\sigma}\Big\|_{E}+
\|(A(w))u\|_{X}.
\end{array}
$$
Since $A(w)$ is the same as the one in the proof of H4 in \cite{AJ},
we have that for each $w\in W$, there exists a constant $c>0$ such
that $ \|(A(w))u\|_{X}\leq c\|u\|_{Y}. $ Then $$
\displaystyle\|\mathbf{A}(\mathbf{W})\mathbf{U}\|_{\mathbb{X}}\leq
(1+c)\|\mathbf{U}\|_{\mathbb{Y}}.
$$ This completes the proof of the first assertion.

For any $\mathbf{W}_{1}, \mathbf{W}_{2}\in
\mathbb{W}$$(\mathbf{W}_{1}=(\tilde{w}_{1},w_{1}),\mathbf{W}_{2}=(\tilde{w}_{2},w_{2}))$,
$$
\|(\mathbf{A}(\mathbf{W}_{1})-\mathbf{A}(\mathbf{W}_{2}))\mathbf{U}\|_{\mathbb{X}}=
\|(A(w_{1})-A(w_{2}))u\|_{X}.
$$
Under (A.1) and (A.3), we can obtain (3.10) by using the same method
in the proof of H4 in \cite{AJ}. This completes the
proof.$\quad\Box$

{\em Remark 3.2:}\ \  Since for a given $\mathbf{W}\in \mathbb{W}$,
$\mathbf{A}(\mathbf{W})$ is an autonomous first-order differential
operator,
$\mathbf{A}(\mathbf{W}):[0,T_{0}]\rightarrow\mathcal{B}(\mathbb{Y},\mathbb{X})$
is strongly continuous with $T_{0}=\infty$.

\medskip
{\bf Lemma 3.5}\ \ {\em There exist positive numbers
$\lambda_{\mathbf{F}}$, $\mu_{\mathbf{F}}$ and
$\bar{\mu}_{\mathbf{F}}$ such that
\begin{eqnarray}
\|\mathbf{F}(\mathbf{W})\|_{\mathbb{Y}}\leq
\lambda_{\mathbf{F}},\;\;\hbox{for}\;\;\mathbf{W}\in \mathbb{W},
\end{eqnarray}
\begin{eqnarray}
\|\mathbf{F}(\mathbf{W}_{1})-\mathbf{F}(\mathbf{W}_{2})\|_{\mathbb{X}}\leq
\mu_{\mathbf{F}}\|\mathbf{W}_{1}-\mathbf{W}_{2}\|_{\mathbb{X}},\;\;\hbox{for}\;\;\mathbf{W}_{1},\mathbf{W}_{2}\in
\mathbb{W},
\end{eqnarray}
\begin{eqnarray}
\|\mathbf{F}(\mathbf{W}_{1})-\mathbf{F}(\mathbf{W}_{2})\|_{\mathbb{Y}}\leq
\bar{\mu}_{\mathbf{F}}\|\mathbf{W}_{1}-\mathbf{W}_{2}\|_{\mathbb{Y}},
\;\;\hbox{for}\;\;\mathbf{W}_{1},\mathbf{W}_{2}\in
\mathbb{W}.
\end{eqnarray} }
\medskip
{\em Proof}:\ \ Let $\mathbf{W}=(\tilde{w},w)\in \mathbb{W}$,
$\mathbf{W}_{1}=(\tilde{w}_{1},w_{1})\in\mathbb{W}$ and
$\mathbf{W}_{2}=(\tilde{w}_{2},w_{2})\in\mathbb{W}$. We have that
\begin{eqnarray}\|\mathbf{F}(\mathbf{W})\|_{\mathbb{Y}}\leq\|f_{1}(w)\|_{Y}+
\|f_{2}(\tilde{w})\|_{Y},
\end{eqnarray}
\begin{eqnarray}\|\mathbf{F}(\mathbf{W}_{1})-\mathbf{F}(\mathbf{W}_{2})\|_{\mathbb{X}}\leq
\|f_{1}(w_{1})-f_{1}(w_{2})\|_{X}+\|f_{2}(\tilde{w}_{1})-f_{2}(\tilde{w}_{1})\|_{X},
\end{eqnarray}
and
\begin{eqnarray}\|\mathbf{F}(\mathbf{W}_{1})-\mathbf{F}(\mathbf{W}_{2})\|_{\mathbb{Y}}\leq
\|f_{1}(w_{1})-f_{1}(w_{2})\|_{Y}+\|f_{2}(\tilde{w}_{1})-f_{2}(\tilde{w}_{1})\|_{Y}.
\end{eqnarray}
Under(A.1) and (A.3), by using the same method in the proof of H5 in
\cite{AJ}, we have that there exist positive numbers
$\lambda_{f_{1}}$, $\mu_{f_{1}}$ and $\bar{\mu}_{f_{1}}$ such that
\begin{eqnarray}
\|f_{1}(w)\|_{Y}\leq \lambda_{f_{1}},
\end{eqnarray}
\begin{eqnarray}
\|f_{1}(w_{1})-f_{1}(w_{2})\|_{X}\leq
\mu_{f_{1}}\|w_{1}-w_{2}\|_{X},
\end{eqnarray}
and
\begin{eqnarray}
\|f_{1}(w_{1})-f_{1}(w_{2})\|_{Y}\leq
\mu_{f_{1}}\|w_{1}-w_{2}\|_{Y}.
\end{eqnarray}
Then we only need to consider $f_{2}$. Under (A.4), we have
\begin{eqnarray}
\displaystyle\|f_{2}(\tilde{w})\|_{Y}=\displaystyle\Big\|\int^{0}_{-\tau}R[\tilde{w}(\sigma,y)]d\sigma
\Big\|_{W^{1,1}(0,\infty)}\leq
R_{1}\int^{\infty}_{0}\int^{0}_{-\tau}\tilde{w}(\sigma,y)d\sigma
dy\leq\mu_{f_{2}},
\end{eqnarray}
\begin{eqnarray}
\displaystyle\|f_{2}(\tilde{w}_{1})-f_{2}(\tilde{w}_{2})\|_{Y}\leq\displaystyle
 L_{R}\|\tilde{w}_{1}-\tilde{w}_{1}\|_{E} \end{eqnarray}
and
\begin{eqnarray}
\displaystyle\|f_{2}(\tilde{w}_{1})-f_{2}(\tilde{w}_{2})\|_{X}\leq\displaystyle
 L_{R_{x}}\|\tilde{w}_{1}-\tilde{w}_{1}\|_{E} \end{eqnarray}

From(3.14-3.22), we obtain (3.11), (3.12) and (3.13). This completes
the proof.$\quad\Box$

{\em Remark 3.3:}\ \ For each $\mathbf{W}\in \mathbb{W}$,
$\mathbf{F}(\mathbf{W})$ is a well defined function belonging to
$\mathbb{Y}$. Since $\mathbf{F}(\mathbf{W})$ does not depend on $t$,
for each $\mathbf{W}\in \mathbb{W}$, $\mathbf{F}(\cdot,\mathbf{W})$
is continuous in $\mathbb{X}$ on $[0,T_{0}]$ and is strongly
measurable in $\mathbb{Y}$ with $T_{0}=\infty$.

By Lemma 3.1-3.5 and Remark 3.1-3.3, the hypotheses of Theorem I and
II in\cite{AJ} are fulfilled with $T_{0} = \infty$ and $\mathbb{W}$
the open subset of $\mathbb{Y}$ which contained in an arbitrary
closed ball in $\mathbb{Y}$ with center $0$ and radius
$r>\|\mathbf{U}_{0}\|_{\mathbb{Y}}$. Hence, we have the following
result:

\medskip
{\bf Lemma 3.6}\ \ {\em Under $(A.1)-(A.4)$, for any initial
condition $\mathbf{U}_{0}\in \mathbb{Y}$, there exists a time $T>0$
such that the problem $(2.5)$ has a unique solution $\mathbf{U}\in
C([0,T],\mathbb{Y})\cap C^{1}([0,T],\mathbb{X})$. Moreover, the
family of operators $\{\mathbb{U}(t,s)\}$, $(t,s)\in \Delta$,
generated by\\ $\{\mathbf{A}(\mathbf{U}(t))\}_{t\in[0,T]}$ is stable
with stability index $(M,\alpha)$ in $\mathbb{X}$ and
$(\tilde{M},\tilde{\alpha})$ in $\mathbb{Y}$, where
$(\tilde{M},\tilde{\alpha})=(M\|\mathbf{S}\|_{\mathbb{Y},\mathbb{X}}
\|\mathbf{S}^{-1}\|_{\mathbb{X},\mathbb{Y}},\lambda_{\mathbf{B}}M+
\alpha)$ and it satisfies the properties of Theorem II in
$\cite{AJ}$.}
\medskip

By Lemma 2.1, we have the following result:

\medskip
{\bf Theorem 3.1}\ \ {\em Under $(A.1)-(A.4)$, for any initial
condition $(\hat{n},\hat{n}_{0})\in W^{1,1}([-\tau,0],X)\times Y$,
there exists a time $T>0$ such that $(2.2)$ has a unique solution
$n\in C([-\tau,T], W^{1,1}\\(0,\infty))\cap C([0,T],Y)\cap
C^{1}([0,T],X)$.}

\section{Continuous Dependence on Initial Conditions and Positivity of Solutions}
\setcounter{equation}{0}

In this section we obtain the continuous dependence on initial
conditions and the positivity of solutions.

Since under $(A.1)-(A.4)$, the hypotheses of Theorem I and II in
\cite{AJ} are fulfilled and the evolution operator of (2.5)
satisfies Theorem II in \cite{AJ}, we obtain the following results
by using the same methods for proving Theorem 2 and Theorem 3 in
\cite{AJ}.

\medskip
{\bf Lemma 4.1}\ \ {\em Let $\mathbf{U}$ and $\mathbf{V}$ be
solutions of the problem (2.5) with initial conditions
$\mathbf{U}_{0}$ and $\mathbf{V}_{0}$, respectively, in
$\mathbb{Y}$. Then, under $(A.1)-(A.4)$, for all $0<t<T$, there
exists a constant $\zeta(r,T)$ such that
$$
\|\mathbf{U}(t)-\mathbf{V}(t)\|_{\mathbb{X}}\leq Me^{\alpha
t}\|\mathbf{U}_{0}-\mathbf{V}_{0}\|_{\mathbb{X}}(1+t\zeta(r,T)),
$$
with
$r>\max\{\|\mathbf{U}_{0}\|_{\mathbb{Y}},\|\mathbf{V}_{0}\|
_{\mathbb{Y}}\}$,
where $T$ is a common local existence time of $\mathbf{U}$ and
$\mathbf{V}$.}
\medskip

\medskip
{\bf Lemma 4.2}\ \ {\em If the initial condition $\mathbf{U}_{0}\geq
0$, then under $(A.1)-(A.4)$, the solution of the problem $(2.5)$ is
non-negative for any $t\in [0,T]$, where $T$ is the local existence
time of the solution.}
\medskip

By Lemma 2.1, we have the following results:

\medskip
{\bf Theorem 4.1}\ \ {\em Let $n$ and $m$ be solutions of $(2.2)$
with initial conditions $(\hat{n},\hat{n}_{0})$ and
$(\hat{m},\hat{m}_{0})$, respectively, in
$W^{1,1}([-\tau,0],X)\times Y$. Then, under $(A.1)-(A.4)$, for all
$0<t<T$, there exists a constant $\zeta(r,T)$ such that
$$
\|n(t)-m(t)\|_{X}\leq Me^{\alpha t}\Big(\|\hat{n}-\hat{m}\|_{E}+
\|\hat{n}_{0}-\hat{m}_{0}\|_{X}\Big)(1+t\zeta(r,T)),
$$
with $r>\max\{\|(\hat{n},\hat{n}_{0})
\|_{\mathbb{Y}},\|(\hat{m},\hat{m}_{0}) \|_{\mathbb{Y}}\}$, where
$T$ is a common local existence time of $n$ and $m$.}
\medskip

\medskip
{\bf Theorem 4.2}\ \ {\em If the initial condition $\hat{n}\geq 0$,
then under $(A.1)-(A.4)$, the solution of the problem $(2.2)$ is
non-negative for any $t\in [-\tau,T]$, where $T$ is the local
existence time of the solution.}
\medskip

\section{Global Existence of Solution}
\setcounter{equation}{0}

We denote by $\varphi(t;t_{0},x_{0})$ the characteristic curve
passing through $(x_{0},t_{0})\in (0,+\infty)\times[0,T]$, i.e., it
is the solution of
\begin{eqnarray}
\frac{\partial \varphi}{\partial
t}=\gamma(\varphi(t;t_{0},x_{0}),N[n(t)](\varphi(t;t_{0},x_{0}))),
\varphi(t_{0};t_{0},x_{0})=x_{0}.
\end{eqnarray}

From (1.1), we have that for $t\geq 0$,
\begin{eqnarray}
n(x,t)&=&\left\{
\begin{array}{l}
\displaystyle \int^{t}_{\eta}
e^{-\int^{t}_{\zeta}\mu(\varphi(s;t,x),N[n(s)](\varphi(s;t,x)))ds}
\\\;\;\;\;\;\;\int^{0}_{-\tau}R[n(\zeta+\sigma)](\varphi(\zeta;t,x))\varphi_{x}(\zeta;t,x)d\zeta, \;\;\;  x<z(t), \\
\displaystyle \hat{n}_{0}(\varphi(0;t,x))e^{-\int^{t}_{0}
\mu(\varphi(s;t,x),N[n(s)](\varphi(s;t,x)))d
s}\varphi_{x}(0;t,x)\\
+\displaystyle
\int^{t}_{0}e^{-\int^{t}_{\zeta}\mu(\varphi(s;t,x),N[n(s)](\varphi(s;t,x)))ds}\\
\;\;\;\;\;\;\int^{0}_{-\tau}R[n(\zeta+\sigma)](\varphi(\zeta;t,x))\varphi_{x}(\zeta;t,x)d\zeta,\;
 x>z(t),
\end{array}
\right.
\end{eqnarray}
where $\eta$ is implicitly given by $\varphi(\eta;t,x)=0$,
$z(t):=\varphi(t;0,0)$ is the characteristic curve coming from the
origin, and
$$\varphi_{x_{0}}(t;t_{0},x_{0}):=\exp\Big(\int^{t}_{t_{0}}
D\gamma(\varphi(s;t_{0},x_{0}),N[n(s)](\varphi(s;t_{0},x_{0})))ds\Big).$$

In order to obtain global existence of solution, we make the
additional assumption :

$(A.5)$ For any positive integrable function $n$,
$\displaystyle\frac{\partial\gamma(x,N)}{\partial N}\frac{\partial
N[n](x)}{\partial x}\geq 0$.

\medskip
{\bf Lemma 5.1}\ \ {\em Let us assume the hypotheses $(A.1)-(A.5)$
and let $\mathbf{U}$ be a positive solution of the problem $(2.5)$
up to time $T$. Then $\|\mathbf{U}\|_{\mathbb{Y}}$ is upper bounded
by a positive continuous function of $t$ for all $t\in [0,T]$.}
\medskip

{\em Proof}:\ \ Since
\begin{eqnarray} \|\mathbf{U}(t)\|_{\mathbb{Y}}=
\displaystyle\|n(t)\|_{X}+\|n_{x}(t)\|_{X}+\|n_{t}\|_{E}+\Big\|\frac{\partial
n(t+\sigma)}{\partial\sigma}\Big\|_{E},\end{eqnarray} we step by
step obtain the estimates of each term.

{\bf The estimate of $\|n(t)\|_{X}$:} From $(5.2)$, we have that
$$
\frac{d\|n(t)\|_{X}}{dt}=\int^{\infty}_{0}
\int^{0}_{-\tau}R[n(t+\sigma,y)](s)d\sigma
-\mu(s,N[n(t)](s))n(s,t))ds.
$$
From the positivity of $n(t)$ and $\mu(x,N)$, we have that
\begin{eqnarray}
\begin{array}{rcl}
\displaystyle\|n(t)\|_{X}&\leq&\displaystyle\int^{t}_{0}\int^{\infty}_{0}
\int^{0}_{-\tau}R[n(r+\sigma)](s)d\sigma d s dr +\|\hat{n}_{0}\|_{X}
\\[0.1cm]&\leq&
\displaystyle
\bar{R}\int^{t}_{0}\int^{\infty}_{0}\int^{0}_{-\tau}n(r+\sigma,y)d\sigma
dy
 dr+\|\hat{n}_{0}\|_{X}\\[0.1cm]&=&
\displaystyle
\bar{R}\int^{\infty}_{0}\int^{0}_{-\tau}\int^{t}_{0}n(r+\sigma,y)dr
d\sigma
 dy+\|\hat{n}_{0}\|_{X}\\[0.1cm]&=&
\displaystyle
\bar{R}\int^{\infty}_{0}\int^{0}_{-\tau}\int^{t+\sigma}_{\sigma}n(\xi,y)d\xi
d\sigma
 dy+\|\hat{n}_{0}\|_{X}\\[0.1cm]
&\leq&\displaystyle\bar{R}\tau\int^{t}_{0}\|n(\xi)\|_{X}d\xi+\bar{R}\tau\|\hat{n}\|_{E}+\|\hat{n}_{0}\|_{X},\nonumber
\end{array}
\end{eqnarray}
where $\bar{R}=\max\{R_{0},R_{1},R_{2}\}$. Using the Gronwall's
lemma, we have that
\begin{eqnarray}
\|n(t)\|_{X}\leq
(\bar{R}\tau\|\hat{n}\|_{E}+\|\hat{n}_{0}\|_{X})e^{\bar{R}\tau t}.
\end{eqnarray}
{\bf The estimate of $\|n_{t}\|_{E}$:} We have that if
$t-\tau\geq 0 $, then
\begin{eqnarray}
\begin{array}{rcl}
\displaystyle\|n_{t}\|_{E}&=&\displaystyle\int^{0}_{-\tau}\int^{\infty}_{0}
n(t+\sigma,y)d\sigma d y
\\[0.1cm]&\leq&
\displaystyle\int^{0}_{-\tau}
(\bar{R}\tau\|\hat{n}\|_{E}+\|\hat{n}_{0}\|_{X})e^{\bar{R}\tau
(t+\sigma)}d\sigma \\[0.1cm]&\leq&\displaystyle \tau(\bar{R}\tau\|\hat{n}\|_{E}+\|\hat{n}_{0}\|_{X})e^{\bar{R}\tau
t},\nonumber
\end{array}
\end{eqnarray}
otherwise, i.e. $t-\tau<0 $, we have
\begin{eqnarray}
\begin{array}{rcl}
\displaystyle\|n_{t}\|_{E}&=&\displaystyle\int^{\infty}_{0}
\int^{-t}_{-\tau}n(t+\sigma,y)d\sigma d y+\int^{\infty}_{0}
\int^{0}_{-t}n(t+\sigma,y)d\sigma d y
\\[0.1cm]&\leq&\displaystyle\int^{\infty}_{0}
\int^{0}_{-\tau}n(\xi,y)d\xi d y+\int^{\infty}_{0}
\int^{t}_{0}n(\xi,y)d\xi d
y\\[0.1cm]&\leq&\displaystyle\|\hat{n}\|_{E}+t(\bar{R}\tau\|\hat{n}\|_{E}+\|\hat{n}_{0}\|_{X})e^{\bar{R}\tau
t}.\nonumber
\end{array}
\end{eqnarray}
Then
\begin{eqnarray}
\displaystyle\|n_{t}\|_{E}\leq\displaystyle
\|\hat{n}\|_{E}+(t+\tau)(\bar{R}\tau\|\hat{n}\|_{E}+\|\hat{n}_{0}\|_{X})e^{\bar{R}\tau
t}.
\end{eqnarray}
{\bf The estimate of $\|n_{x}(t)\|_{X}$:} Let
$\tilde{n}(t):=n(\varphi(t;t_{0},x_{0}),t)$. We first give the
estimate of $\|\tilde{n}(t)\|_{\infty}$. Let
$\tilde{n}(t):=n(\varphi(t;t_{0},x_{0}),t)$. We consider $(1.1)_{1}$
as a directional derivative and, so for a fixed $(x_{0},t_{0})$, we
have that
$$
\frac{d\tilde{n}(t)}{dt}=
\int^{0}_{-\tau}R[n(t+\sigma)]\varphi(t;t_{0},x_{0})d\sigma
-\lambda(\varphi(t;t_{0},x_{0}),t)\tilde{n}(t)dt+\tilde{n}(0),
$$
where $\lambda(x,t)=D\gamma(x,N[n(t)(x)])+\mu(s,N[n(t)](x))$. Then
$$
\tilde{n}(t)=\frac{d\tilde{n}(t)}{dt}=\int^{t}_{0}
\int^{0}_{-\tau}R[n(r+\sigma)]\varphi(r;t_{0},x_{0})d\sigma dr
-\int^{t}_{0}\lambda(\varphi(r;t_{0},x_{0}),t)\tilde{n}(r)dr+\tilde{n}(0),
$$
Recall that $\displaystyle\frac{\partial\gamma(x,N)}{\partial
N}\frac{\partial N[n](x)}{\partial x}$ and $\mu(x,N)$ are
non-negative. We denote
$\lambda_{0}=inf_{x,t}\\\{\gamma_{x}(x,N(x,t)])\}$ and assume that
$\lambda_{0}$ is negative, and $\lambda_{0}=0$ otherwise. Then
\begin{eqnarray}
\begin{array}{rcl}
\displaystyle\tilde{n}(t)&\leq&
\displaystyle\bar{R}\int^{t}_{0}\int^{\infty}_{0}\int^{0}_{-\tau}n(r+\sigma,y)d\sigma
dydr-\lambda_{0}\int^{t}_{0}\tilde{n}(r)dr+\tilde{n}(0)\\[0.1cm]&\leq&
\displaystyle\bar{R}\int^{\infty}_{0}\int^{0}_{-\tau}\int^{t+\sigma}_{\sigma}n(\xi,y)d\xi
d\sigma dy-\lambda_{0}\int^{t}_{0}\tilde{n}(r)dr+\tilde{n}(0)\\[0.1cm]&\leq&
\displaystyle
\bar{R}\tau\|\hat{n}\|_{E}+\bar{R}\tau\int^{t}_{0}\|n(\xi)\|_{X}d\xi
-\lambda_{0}\int^{t}_{0}\tilde{n}(r)dr+\tilde{n}(0)\\[0.1cm]&\leq&
\displaystyle\bar{R}\tau\|\hat{n}\|_{E}+\tilde{n}(0)+
(\bar{R}\tau\|\hat{n}\|_{E}+\|\hat{n}_{0}\|_{X})e^{\bar{R}\tau t}
-\lambda_{0}\int^{t}_{0}\tilde{n}(r)dr.\nonumber
\end{array}
\end{eqnarray}
Using the Gronwall's lemma, we have that if $\lambda_{0}\neq
-\bar{R}\tau $, then
\begin{eqnarray}
\displaystyle\|\tilde{n}(t)\|_{\infty}\leq
\displaystyle(\bar{R}\tau\|\hat{n}\|_{E}+\|\hat{n}_{0}\|_{\infty})(e^{-\lambda_{0}t}+1)+
(\bar{R}\tau\|\hat{n}\|_{E}+\|\hat{n}_{0}\|_{X}) e^{\bar{R}\tau
t}(1-\frac{\lambda_{0}}{\bar{R}\tau+\lambda_{0}})=:g_{1}(t),\nonumber\\
\end{eqnarray}
otherwise, i.e. $\lambda_{0}=-\bar{R}\tau<0 $, we have
\begin{eqnarray}
\displaystyle\|\tilde{n}(t)\|_{\infty}\leq
\displaystyle(\bar{R}\tau\|\hat{n}\|_{E}+\|\hat{n}_{0}\|_{\infty})
(e^{-\lambda_{0}t}+1)+(\bar{R}\tau\|\hat{n}\|_{E}+\|hat{n}_{0}\|_{X})
(e^{\bar{R}\tau t}-\lambda_{0}te^{ -\lambda_{0}t})=:g_{2}(t).\nonumber\\
\end{eqnarray}
Let us consider $n$ as a function of $t$ and $\xi:=\varphi(0;t,x)$
for $x>z(t)$, and as a function of $t$ and $\eta$, with $\eta$ given
by $\varphi(\eta;t,x)=0$ for $x<z(t)$, i.e.
$$
n(x,t)=\bar{n}(t,\eta(x,t)),0<x<z(t),n(x,t)=\bar{\bar{n}}(t,\xi(x,t)),x>z(t).
$$
Hence, we have that
\begin{eqnarray}
n_{x}(x,t)&=&\left\{
\begin{array}{l}
\displaystyle\bar{n}_{\eta}(t,\eta)\partial_{x}\eta(x,t), \;\;\;  x<z(t), \\
\displaystyle \bar{\bar{n}}_{\xi}(t,\xi)\partial_{x}\xi(x,t),\;
 x>z(t),
\end{array}
\right.\nonumber
\end{eqnarray}
 and
\begin{eqnarray}
\begin{array}{rcl}
\displaystyle\|n_{x}(t)\|_{X}&=&\displaystyle\int^{z(t)}_{0}
|n_{x}(x,t)|dx+\int^{\infty}_{z(t)} |n_{x}(x,t)|dx
\\[0.1cm]&\leq&
\displaystyle \int^{t}_{0}
|\bar{n}_{\eta}(t,\eta)|d\eta+\int^{\infty}_{0}
|\bar{\bar{n}}_{\xi}(t,\xi)|d\xi.\nonumber
\end{array}
\end{eqnarray}
The expression of $\bar{n}_{\eta}$ and $\bar{\bar{n}}_{\xi}$ are
similar as $\bar{u}_{\tau}$ and $\bar{\bar{u}}_{\xi}$ in the proof
of Lemma 1 in \cite{AJ}, we obtain the following result by using the
same method:
\begin{eqnarray}
\displaystyle \|n_{x}(t)\|_{X}\leq \displaystyle
f_{1}(t)+\displaystyle[\|\hat{n}_{0}\|_{\infty}e^{-\lambda_{0}t}+
g(t)](\gamma^{0}_{2}\int^{t}_{0}\|(N_{x}(s))^{2}\|_{X}d
s+\gamma^{0}_{1}\int^{t}_{0}\|N_{xx}(s)\|_{X}d s).\nonumber\\
\end{eqnarray}
where $f_{1}(t)$ is a positive and increasing function which is
similar as the one in (6.8) of \cite{AJ}, $g(t):=g_{1}(t)+g_{2}(t)$.
Moreover, under (A.1) and (A.2) and using the $\|N_{x}\|_{L^{1}}\leq
c^{'}\|n(t)\|_{L^{1}}$ and
$\|N_{x}\|_{\infty}\leq\|N_{x}\|_{W^{1,1}}\leq
c^{''}\|n(t)\|_{W^{1,1}}$, the last two integrals of (5.8) can be
bounded as follows
\begin{eqnarray}
\begin{array}{rcl}
\displaystyle \int^{t}_{0}\|(N_{x}(s))^{2}\|_{X}d
s&\leq&\displaystyle
\int^{t}_{0}\|N_{x}(s)\|_{\infty}\|N_{x}(s)\|_{X}ds
\\[0.1cm]&\leq&\displaystyle c'(\bar{R}\tau\|\hat{n}\|_{E}+
\|\hat{n}_{0}\|_{X})e^{\bar{R}\tau t}
\int^{t}_{0}\|N_{x}(s)\|_{\infty}ds
\\[0.1cm]&\leq&
\displaystyle
c'(\bar{R}\tau\|\hat{n}\|_{E}+\|\hat{n}_{0}\|_{X})e^{\bar{R}\tau
t}c''\int^{t}_{0}(\|n(s)\|_{X}+\|n_{x}(s)\|_{X})ds
\end{array}
\end{eqnarray}
and \begin{eqnarray} \displaystyle \int^{t}_{0}\|N_{xx}(s)\|_{X}d
s\leq \displaystyle
c''\int^{t}_{0}(\|n(s)\|_{X}+\|n_{x}(s)\|_{X})ds.
\end{eqnarray}
From (5.8)-(5.10), we have that
\begin{eqnarray}
\begin{array}{rcl}
\displaystyle \|n_{x}(t)\|_{X}&\leq&\displaystyle
f(t)+f_{2}(t)\int^{t}_{0}\|n(r)\|_{X} dr,
\end{array}
\end{eqnarray}
where $f_{2}(t)=c''[\|\hat{n}_{0}\|_{\infty}e^{-\lambda_{0}t}+g(t)]
(\gamma^{0}_{2}c'(\bar{R}\tau\|\hat{n}\|_{E}+\|\hat{n}_{0}\|_{X})e^{\bar{R}\tau
t}+\gamma^{0}_{1})$ and
$f(t)=f_{1}(t)+tf_{2}(t)\bar{R}\tau\|\hat{n}\|_{E}+\|\hat{n}_{0}\|_{X})e^{\bar{R}\tau
t}$.

From the Gronwall's lemma, we have that
\begin{eqnarray}
\begin{array}{rcl}
\displaystyle \|n_{x}(t)\|_{X}&\leq&\displaystyle
f(t)+f_{2}(t)\int^{t}_{0}f(s)\exp\Big(\int^{t}_{s}f_{2}(\zeta)d\zeta\Big)ds.
\end{array}
\end{eqnarray}
Then, we have that $\|n_{x}(t)\|_{X}\leq
H_{1}(t,\|\mathbf{U}_{0}\|_{\mathbb{Y}})$, where
$H_{1}(t,\|\mathbf{U}_{0}\|_{\mathbb{Y}})$ denotes the R.H.S of
(5.13) with the $\|\hat{n}\|_{E}$, $\|\hat{n}_{0}\|_{\infty}$ and
$\|\hat{n}_{0}\|_{L^{1}}$ replaced by
$\|\mathbf{U}_{0}\|_{\mathbb{Y}}$.

{\bf The estimate of $\displaystyle\Big\|\frac{\partial
n(t+\sigma)}{\partial\sigma}\Big\|_{E}$:}
\begin{eqnarray}
\displaystyle\int^{\infty}_{0}\int^{0}_{-\tau}\Big|\frac{\partial
n(t+\sigma,x)}{\partial\sigma}\Big|d\sigma
dx=\int^{\infty}_{0}\int^{t}_{-\tau+t}\Big|\frac{\partial
n(\xi,x)}{\partial\xi}\Big|d\xi dx.\nonumber
\end{eqnarray}
If $t>\tau$, $n(\xi,x)$ satisfies $(1.1)_{1}$, then we have that
\begin{eqnarray}
&& \displaystyle
\int^{\infty}_{0}\int^{t}_{-\tau+t}\Big|\frac{\partial
n(\xi,x)}{\partial\xi}\Big|d\xi dx\nonumber
\\&\leq&\displaystyle
\int^{\infty}_{0}\int^{t}_{0}\Big|\frac{\partial
n(\xi,x)}{\partial\xi}\Big|d\xi dx\nonumber
\\&\leq&
\displaystyle\gamma^{0}\int^{t}_{0}\int^{\infty}_{0}\Big|\frac{\partial
n(\xi,x)}{\partial x}\Big|dx
d\xi\displaystyle+(\mu^{0}+\gamma^{0}_{1})\int^{t}_{0}
\int^{\infty}_{0}|n(\xi,x)|dx
d\xi\nonumber\\[0.1cm]&&+ \displaystyle
\int^{t}_{0}\int^{\infty}_{0}\int^{0}_{-\tau}
R[n(\xi+\sigma)](x)d\sigma dx d\xi\nonumber\\[0.1cm]&\leq&
\displaystyle\gamma^{0}t
N_{2}(t,\|\mathbf{U}_{0}\|_{\mathbb{Y}})+(\mu^{0}+\gamma^{0}_{1})t(\bar{R}\tau\|\hat{n}\|_{E}+
\|\hat{n}_{0}\|_{X})e^{\bar{R}\tau t})+
\bar{R}\int^{t}_{0}\int^{\infty}_{0}\int^{0}_{-\tau}n(\xi+\sigma,y)d\sigma
dyd\xi
\nonumber\\[0.1cm]&\leq&
\displaystyle\gamma^{0}t H_{1}(t,\|\mathbf{U}_{0}\|_{\mathbb{Y}})+
(\mu^{0}+\gamma^{0}_{1})t(\bar{R}\tau\|\hat{n}\|_{E}+
\|\hat{n}_{0}\|_{X})e^{\bar{R}\tau t})\nonumber\\[0.1cm]&&+
\bar{R}t(\|\hat{n}\|_{E}+(t+\tau)(\bar{R}\tau\|\hat{n}\|_{E}+\|\hat{n}_{0}\|_{X})e^{\bar{R}\tau
t}),\nonumber
\end{eqnarray}
otherwise, i.e. $t\leq\tau$, we have
\begin{eqnarray}
&& \displaystyle
\int^{\infty}_{0}\int^{t}_{-\tau+t}\Big|\frac{\partial
n(\xi,x)}{\partial\xi}\Big|d\xi dx\nonumber
\\&\leq&\displaystyle\int^{\infty}_{0}\int^{0}_{-\tau+t}\Big|\frac{\partial
n(\xi,x)}{\partial\xi}\Big|d\xi dx+
\int^{\infty}_{0}\int^{t}_{0}\Big|\frac{\partial
n(\xi,x)}{\partial\xi}\Big|d\xi dx\nonumber
\\&\leq&
\displaystyle\int^{\infty}_{0}\int^{0}_{-\tau}\Big|\frac{\partial
\hat{n}(\xi,x)}{\partial\xi}\Big|d\xi dx+
\int^{\infty}_{0}\int^{t}_{0}\Big|\frac{\partial
n(\xi,x)}{\partial\xi}\Big|d\xi dx
\nonumber\\[0.1cm]&\leq&
\displaystyle\|\mathbf{U}_{0}\|_{\mathbb{Y}}+\gamma^{0}t
H_{1}(t,\|\mathbf{U}_{0}\|_{\mathbb{Y}})+(\mu^{0}+\gamma^{0}_{1})t(\bar{R}\tau\|\hat{n}\|_{E}+\|\hat{n}_{0}\|_{X})e^{\bar{R}\tau
t})\nonumber\\[0.1cm]&&+
\bar{R}t(\|\hat{n}\|_{E}+(t+\tau)(\bar{R}\tau\|\hat{n}\|_{E}+\|\hat{n}_{0}\|_{X})e^{\bar{R}\tau
t}).\nonumber
\end{eqnarray}
Then
\begin{eqnarray}
\begin{array}{rcl}
\displaystyle\Big\|\frac{\partial
n(t+\sigma)}{\partial\sigma}\Big\|_{E}&\leq&\displaystyle
\|\mathbf{U}_{0}\|_{\mathbb{Y}}+\gamma^{0}t
H_{1}(t,\|\mathbf{U}_{0}\|_{\mathbb{Y}})+(\mu^{0}+\gamma^{0}_{1})t(\bar{R}\tau\|\hat{n}\|_{E}+\|\hat{n}_{0}\|_{X})e^{\bar{R}\tau
t})\\[0.1cm]&&+
\bar{R}t(\|\hat{n}\|_{E}+(t+\tau)(\bar{R}\tau\|\hat{n}\|_{E}+\|\hat{n}_{0}\|_{X})e^{\bar{R}\tau
t}).
\end{array}
\end{eqnarray}
From (5.3), (5.4), (5.5), (5.12) and (5.13), we have the proof of
this theorem. $\quad\Box$

Using the same method in the Section 6.3 of \cite{AJ}, we have the
following result:

\medskip
{\bf Lemma 5.2}\ \ {\em Under $(A.1)-(A.5)$ and for any
$T^{^{\ast}}>0$, the problem $(2.5)$ has a unique solution up to
$T^{^{\ast}}$, which is positive whenever the initial condition
$\mathbf{U}_{0}>0$.}
\medskip

By Lemma 2.1, we have the following result:

\medskip
{\bf Theorem 5.3}\ \ {\em Under $(A.1)-(A.5)$ and for any
$T^{^{\ast}}>0$, the problem $(2.2)$ has a unique solution up to
$T^{^{\ast}}$, which is positive whenever the initial condition
$\hat{n}>0$.}
\medskip

\section{Remark}\
In this section, we give the typical examples of the operators $N$
and $R$. The environment experienced by an individual of size $x$
when the population density is $n(x,t)$ can be given by
$$
N[n(t)](x)=\int^{\infty}_{0}\rho(x,y)n(t,y)dy.
$$
Many examples of such a sort of environments with different $\rho$
can be seen, for instance, in the references \cite{AJ}, \cite{AKS}
\cite{KWB}, \cite{AA}, \cite{AAAD}, \cite{MJ} and\cite{EAK}. A
simple example of the operator $R$ can be given by
$$
R[n(t+\sigma)](x)=\int^{\infty}_{0}\beta(\sigma,x,y)n(t+\sigma,y)d\sigma
dy
$$
where $\beta(\sigma,x,y)$ denotes the reproduction rate which
individuals of size $y$ give birth to the individuals of size $x$
after a time lag $-\sigma$ starting from conception.
$\beta(\sigma,x,y)$ is supposed to satisfy the following conditions:

$(H.1)$ $\beta\in C([-\tau,0]\times[0,\infty)\times[0,\infty))$,
$\beta\geq0$ and is uniformly bounded by $R_{2}$.

$(H.2)$ For all $\sigma,y\in [-\tau,0]\times[0,\infty)$,
$\int^{\infty}_{0}\beta(\sigma,x,y)dx\leq R_{0}$,
$\beta(\sigma,\cdot,y)\in C^{1}[0,\infty)$ and
$\int^{\infty}_{0}\beta_{x}(\sigma,x,y)dx\leq R_{1}$.

In this paper we consider the case in which the process of the
recruitment is influenced by the environment from start to finish,
hence we assume that the rate $\beta$ is also dependent on the total
environment in the time lag between the beginning and the end of the
recruitment. Under such an assumption, the operator $R$ can be given
by
$$
R[n(t+\sigma)](x)=\int^{\infty}_{0}\beta(\sigma,x,y,\mathcal{N}
(\sigma,x))n(t+\sigma,y)d\sigma dy,
$$
where $$
\mathcal{N}(\sigma,x)=\int^{0}_{-\tau}\chi(\sigma,\xi)N[n(t+\xi)](x)d\xi
$$
and
\begin{eqnarray}
\chi(\sigma,\xi)&=&\left\{
\begin{array}{l}
\displaystyle 1, \;\;\;  \xi\geq \sigma, \\
\displaystyle 0, \;\;\;  \xi<\sigma.
\end{array}
\right.\nonumber
\end{eqnarray}
 $\beta(\sigma,x,y, \mathcal{N})$ is supposed to satisfy the
following conditions:

$(H.3)$ $\beta\in
C([-\tau,0]\times[0,\infty)\times[0,\infty)\times[0,\infty))$,
$\beta\geq0$ and is uniformly bounded by $R_{2}$.

$(H.4)$ For all $\sigma,y,\mathcal{N}\in
[-\tau,0]\times[0,\infty)\times[0,\infty)$,
$\int^{\infty}_{0}\beta(\sigma,x,y,\mathcal{N})dx\leq R_{0}$,
$\beta(\sigma,\cdot,y,\mathcal{N})\in C^{1}[0,\infty)$ and
$\int^{\infty}_{0}\beta_{x}(\sigma,x,y,\mathcal{N})dx\leq R_{1}$.
Moreover, $\int^{\infty}_{0}\beta(\sigma,x,y,\mathcal{N})dx$ and
$\int^{\infty}_{0}\beta_{x}(\sigma,x,y,\mathcal{N})dx$ are
Lipschizian functions with respect to $\mathcal{N}$.

{\small

\end{document}